\documentclass[12pt]{article}
\usepackage{verbatim}
\usepackage[T1]{fontenc}
\usepackage[utf8]{inputenc}

\usepackage{amsmath,amsfonts,amssymb,amscd,mathrsfs,latexsym,epsfig}

\usepackage{eufrak}
\newcommand{\C}{{\EuFrak C}}

\newcommand{\bi}{\begin{itemize}}
\newcommand{\ei}{\end{itemize}}

\def\Ind#1#2{#1\setbox0=\hbox{$#1x$}\kern\wd0\hbox to 0pt{\hss$#1\mid$\hss}
\lower.9\ht0\hbox to 0pt{\hss$#1\smile$\hss}\kern\wd0}

\def\Notind#1#2{#1\setbox0=\hbox{$#1x$}\kern\wd0\hbox to 0pt{\mathchardef
\nn=12854\hss$#1\nn$\kern1.4\wd0\hss}\hbox to
0pt{\hss$#1\mid$\hss}\lower.9\ht0 \hbox to
0pt{\hss$#1\smile$\hss}\kern\wd0}

\newcommand{\acl}{\operatorname{acl}}
\newcommand{\tp}{\operatorname{tp}}

\newcommand{\id}{\operatorname{id}}
\newcommand{\LIM}{\operatorname{LIM}}
\newcommand{\Th}{\operatorname{Th}}
\newcommand{\alg}{\operatorname{alg}}

\newtheorem{theorem}{Theorem}[section]
\newtheorem{lemma}[theorem]{Lemma}
\newtheorem{fact}[theorem]{Fact}

\newtheorem{corollary}[theorem]{Corollary}
\newtheorem{proposition}[theorem]{Proposition}
\newtheorem{definition}[theorem]{Definition}
\newtheorem{remark}[theorem]{Remark}

\newtheorem{question}[theorem]{Question}

\newtheorem{main theorem}{Theorem}
\newtheorem{main lemma}[main theorem]{Lemma}
\newtheorem{main fact}[main theorem]{Fact}
\newtheorem{main corollary}[main theorem]{Corollary}
\newtheorem{main proposition}[main theorem]{Proposition}
\newtheorem{main definition}[main theorem]{Definition}
\newtheorem{main remark}[main theorem]{Remark}
\newtheorem{main conjecture}[main theorem]{Conjecture}
\newtheorem{main question}[main theorem]{Question}
\newtheorem{main problem}[main theorem]{Problem}

\newcounter{licznik}
\newtheorem{proven main theorem}[licznik]{Theorem}

\title{On regular groups and fields}
\author{Tomasz Gogacz and Krzysztof Krupi\'nski\footnote{Research supported by the Polish Government grant N N201 545938}}
\date{}

\begin{document}
\maketitle
\begin{abstract}
Regular groups and fields are common generalizations of minimal and quasi-minimal groups and fields, so the conjectures that minimal or quasi-minimal fields are algebraically closed have their common generalization to the conjecture that each regular field is algebraically closed. 
Standard arguments show that a generically stable regular field is algebraically closed. Let $K$ be a regular field which is not generically stable and let $p$ be its global generic type. We observe that if $K$ has a finite  extension $L$ of degree $n$, then $p^{(n)}$ has unbounded orbit under the action of the multiplicative group of $L$.
%We make some observations on regular fields which are not generically stable. If such a filed $K$ has a finite  extension $L$ of degree $n$, then: $p^n$ has unbounded orbit under the multiplicative action of $L$, where $p$ is the unique global generic type of $K$; assuming NIP,  if $q$ is a $\mu$-generic type in $n$-variables with respect to a sufficiently invariant measure $\mu$, then $q$ is invariant under permuting variables and has bounded orbit under the multiplicative action of $L$.

Known to be true in the minimal context, it remains wide open whether regular, or even quasi-minimal,  groups are abelian. We show that if it is not the case, then there is a counter-example with a unique non-trivial conjugacy class, and we notice that a classical group with one non-trivial conjugacy class is not quasi-minimal, because the centralizers of all elements are uncountable. Then we construct a group of cardinality $\omega_1$ with only one non-trivial conjugacy class and such that the centralizers of all non-trivial elements are countable.
\end{abstract}
\footnotetext{2010 MSC: Primary 03C60, Secondary 12L12, 20A15, 20E06, 03C45}
\footnotetext{Key words and phrases: minimal, quasi-minimal, group, field, regular type, generically stable type, HNN-extension}

\section{Introduction}

Recall that a minimal structure is an infinite structure whose all definable subsets are finite or co-finite. A quasi-minimal structure is an uncountable structure in a countable language whose all definable subsets are countable or co-countable. Among the main examples of quasi-minimal structures there are various orders (e.g. $\omega_1 \times \mathbb Q$), strongly minimal structures expanded by some orders, and Zilber's pseudoexponential fields \cite{Z}. A well-known conjecture of Boris Zilber predicts that the complex exponential field $(\mathbb C, +, \cdot,0,1, \exp)$ is quasi-minimal.

A fundamental theorem of Reineke \cite{R} tells us that each minimal group is abelian. Surprisingly, an analogous statement for quasi-minimal groups seems hard to prove.

\begin{main conjecture}\label{conjecture:1}
Each quasi-minimal  group is abelian.
\end{main conjecture} 

A more general question has been formulated in \cite[Section 3]{PT} (see Section 1 for the definition of regularity).

\begin{main question}\label{question:2}
Is every regular group abelian?
\end{main question}

In Section 3 of this paper, we reduce the problem to the case of groups with only one non-trivial conjugacy class. Then we notice that a standard construction (involving HNN-extensions) of an uncountable group with a unique non-trivial conjugacy class  does not lead to a quasi-minimal group, because the centralizers of all non-trivial elements of the resulting group are uncountable (and so also co-uncountable). Motivated by this obstacle, we construct a group of cardinality $\omega_1$ with only one non-trivial conjugacy class, in which the centralizers of all non-trivial elements are countable. We leave as  an open question whether the group we constructed is quasi-minimal, or at least regular. 

One of the oldest unsolved problems in algebraic model theory is Podewski's conjecture predicting that each minimal field is algebraically closed. Known to be true in positive characteristic \cite{W}, it remains wide open in the zero characteristic case. It was recently noticed in \cite{KTW} that it holds when the generic type of the field in question is generically stable (see Section 1 for definitions); in the non generically stable situation, some partial results have been obtained (e.g. elimination of the so-called almost liner case). 

One has an obvious analog of Podewski's conjecture for quasi-minimal fields.

\begin{main conjecture}\label{conjecture:quasi-minimal fields}
Each quasi-minimal field is algebraically closed.
\end{main conjecture} 

The above conjecture is open even in positive characteristic. A common generalization of Podewski's conjecture and Conjecture \ref{conjecture:quasi-minimal fields} is

\begin{main conjecture}\label{conjecture:regular fields}
Each regular field is algebraically closed.
\end{main conjecture}

Applying either the proof of \cite[Theorem 1]{KTW} or of \cite[Theorem 1.13]{HLS}, one easily gets

\begin{main theorem}\label{main theorem: gen stab}
Each generically stable regular field is algebraically closed. In particular, each generically stable minimal or quasi-minimal field is algebriacally closed.
\end{main theorem}

As a consequence, one gets that each quasi-minimal field of cardinality greater than $\omega_1$ is algebraically closed, and a similar result for regular fields with NSOP. The case of minimal or quasi-minimal fields with NIP is still open.

Having Theorem \ref{main theorem: gen stab}, a natural question arises whether there exists a regular field which is not generically stable. It is conjectured in \cite{KTW} that there are no such minimal fields, but in the quasi-minimal context such a filed exists \cite[Example 5.1]{PT}, and, moreover, it is almost linear using the terminology from \cite{KTW}. On the other hand,  it was proved in \cite{KTW} that there are no almost linear minimal fields.  This reveals a certain difference between minimal and quasi-minimal fields.

In the main result of Section 2, we study the situation when the filed in question is not genericllay stable.
So assume now that $K$ is a regular field which is not generically stable. Suppose that $K$ is not algebrically closed, i.e., it has a finite extension $L$ of degree $n$. Then $L$ is naturally interpreted as $K^n$ with coordinate-wise addition and some definable multiplication. Let $p$ be the global generic type of $K$ and $p^{(n)}$ its $n$-th power. We prove that the orbit of $p^{(n)}$ under the multiplicaive group of $L$ is unbounded, which one can hope to be useful to get a final contradiction for some (e.g. NIP) fields. 
More detailed discussion concerning this appears in Section 2.
%More detailed motivation and some other related results are given in Section 2.

The topic of this paper was undertaken  in 2009, when Predrag Tanovi\'c and the second author made independently various observations and found examples concerning minimal and quasi-minimal groups and fields, in particular the first example of a non generically stable quasi-minimal group (modified later to the context of fields in \cite{PT}). Some of these results are contained in \cite{PT} and \cite{KTW}, some other ideas were further developed and contained in the first author's Master's thesis completed in June 2011. The current paper is based on this thesis.

The second author would like to thank Clifton Ealy for sharing useful ideas concerning quasi-minimal groups. The observation that a classical group with only one non-trivial conjugacy class is not quasi-minimal was first made by Clifton Ealy.
%In particular, it was first noticed by him that a classical group with only one non-trivial conjugacy class is not quasiminial. 

\section{Definitions and basic facts}

Throughout this section, $\C$ denotes a monster model of a first order theory $T$,  i.e., a $\kappa$-saturated and  $\kappa$-strongly homogeneous model for large enough cardinal $\kappa$. 
A global type is a type over $\C$. We say that the set is bounded (or small) if it has cardinality smaller than $\kappa$.

The following notion of strongly regular types in an arbitrary theory has been introduced in \cite{PT}; it coincides with the notion of strongly regular types in the stable context.

\begin{definition}
Let $p(x)$ be a global non-algebraic type. Suppose $\varphi(x) \in p(x)$. We say that $(p(x), \varphi(x))$ is strongly regular if, for some
small $A$ over which $\varphi$ is defined, $p$ is $A$-invariant and for all $B\supseteq A$ and $a$
satisfying $\varphi(x)$, either $a \models p |_B$ or $p |_B \vdash p |_{B,a}$.
\end{definition}

If $(p(x), \varphi(x))$ is strongly regular, then as a witness-set $A$ one can take any set over which $p$ is invariant and $\varphi$ is defined.

\begin{definition}
Let $p$ be a global type invariant over $A$.
 A Morley sequence in $p$ over $A$ is a sequence $\langle a_\alpha\rangle_\alpha$ such that $ a_\alpha \models  p|_{Aa_{<\alpha}}$ for all $\alpha$. The invariant defining scheme for $p$ uniquely determines extensions of $p$ to global types over bigger monster models, so we can talk about a Morley sequence in $p$ over $\C$ which will be just called a Morley sequence in $p$. By $ p^{(n)}$ we denote $\tp (b_0, \ldots, b_{n-1}/\C)$, where $(b_0, \ldots, b_{n-1})$ is Morley sequence in $p$.
\end{definition}

If $p$ is an $A$-invariant type, then all Morley sequences of a fixed length in $p$ over $A$ have the same type over $A$ (so the definition of $p^{(n)}$ is correct), and they are indiscernible over $A$.

\begin{definition}
A type $p \in S(\C)$ is generically stable if, for a small set $A$, $p$  is $A$-invariant 
and for every formula $\phi(x)$ (with parameters from $\C$) there exists a natural number $N$ such that for every Morley sequence $(a_i)_{i<\omega}$ in $p$ over $A$
the set  $\{i \in \omega : \models \phi(a_i)\}$ or its complement has cardinality less than $N$.
\end{definition}

\begin{definition}
Let $M$ be any submodel of $\C$ (possibly $M = \C$), and let
$p \in S_1(M)$. The operator $cl_p$ is defined on (all) subsets of $M$ by:
$cl_p(X) = \{a \in M : a \not\models p | X\}$.
Also, define $cl^A_p
(X) = cl_p(X \cup A)$ for any $A \subset M$ and $X \subseteq M$.
\end{definition}

By \cite[Lemma 3.1, Corollary 3.1]{PT}, one has

\begin{fact}
Let $p\in S_1(\C)$ be a non-algebraic type invariant over $A$.\\
(i) $(p(x), x = x)$ is strongly regular iff $cl^A_p$
is a closure operator on $\C$.\\
(ii) Suppose that $(p(x), x = x)$ is strongly regular. Then $cl^A_p$
is a pregeometry
operator on $\C$ iff every Morley sequence in $p$ over $A$ is totally indiscernible over $A$ iff every Morley sequence in $p$ (over $\C$) is totaly indiscernible over $\C$ iff $p$ is generically stable.
\end{fact}

\cite[Theorem 3.1]{PT} yields the following dichotomy for a strongly regular type $p\in S_1(\C)$ invariant over $A$: either $p$ is generically stable, or for some $A_0 \supseteq A$ there is a definable partial order on $\C$ such that every Morley sequence in $p$ over $A_0$ is strictly increasing. In particular, we have

\begin{corollary}\label{corollary:SOP}
If a strongly regular type is not generically stable, then the theory $T$ has SOP.
\end{corollary}

Now, we recall the notion of a regular group.

\begin{definition}
Let $G$ be a definable group in $\C$. Then $G$ is called a regular group
if for some global type $p(x) \in S_G(\C)$, $(p(x), ``x \in G")$ is strongly regular (in
particular, invariant over some small set). A field $K$ definable in $\C$ is regular if its additive (equivalently multiplicative) group is regular.
\end{definition}

In \cite{PT}, also local strong regularity is considered. However, since it is proved in \cite[Theorem 7.2]{PT} that a ``locally regular group'' is regular, we will not talk about local regularity in this paper.
The following fact is \cite[Theorem 3.2]{PT}.

\begin{theorem}
Suppose that $G$ is a group definable over $\emptyset$, which is regular,
witnessed by $p(x) \in S_G(\C)$. Then:\\
(i) $p(x)$ is both left and right translation invariant (and, in fact, invariant
under definable bijections).\\
(ii) a formula $\varphi(x)$ is in $p(x)$ iff two left [right] translates of $\varphi(x)$ cover $G$ iff
finitely many left [right] translates of $\varphi(x)$ cover $G$. Hence, $p(x)$ is the
unique (global) generic (in the sense of translates) type of $G$.\\
(iii) $p(x)$ is definable over $\emptyset$.\\
(iv) $G = G^0$ (i.e., G is connected).
\end{theorem}

We say that an element $a$ from a regular group $G$ is generic over $A\subset \C$ if $a\models p | A$, where $p$ is the unique generic type of $G$. 

\begin{corollary}\label{corollary:dod}
Let $K$ be a regular field with prime field $F$. If $a \in K$ is a generic over $A\subset K$, elements $c,d$ belong to $F(A)$ and $c \neq 0$, then $ac+d$ is generic over $A$.
\end{corollary}

\begin{definition}
A regular group is said to be generically stable, if its unique global generic type is generically stable (over $\emptyset$), and similarly for fields.
\end{definition}

Assume $M$ is a minimal [or quasi-minimal] structure. Let $p \in S_1(M)$ be the type consisting of all formulas over $M$ which define co-finite [co-countable, respectively] subsets of $M$. Not accidentally, this type is called the generic type of $M$. Namely, assume now that $(M, \cdot,\dots)$ is a minimal or quasi-minimal group and $\varphi(x)$ is any formula;  then $\varphi(x) \in p$ iff $\varphi(M) \cdot \varphi(M)=M$ iff two left 
 [right] translates of $\varphi(M)$ cover $M$ iff finitely many left [right] trnaslates of $\varphi(M)$ cover $M$. In particular, $p$ is definable over $\emptyset$, and so it has a unique global heir $\hat p$ in $S(\C)$ (where $\C \succ M$ is a monster model). It is easy to check that $cl_p$ is a closure operator on $M$, so $cl_{\hat p}$ is a closure operator on $\C$ by the definability of $p$, and thus $(\hat p, x=x)$ is strongly regular (see \cite[Theorem 5.3]{PT}).

\begin{corollary}
A monster model of a minimal or quasi-minimal group is regular.
\end{corollary}

Finally, we recall fundamental issues concerning HNN-extensions, which will be used in the last section. For more details the reader is referred to \cite[Chapter IV]{LS}.

\begin{fact}\label{pf}
If G is a group, $A,B \leq G $ and $\varphi: A \rightarrow B$ is an isomorphism, then the group $F$ defined by the presentation
$$
\langle G, t | \forall a \in A (tat^{-1}= \varphi(a))\rangle
$$ 
contains $G$ as a subgroup, and it is called the HNN-extension of $G$ relative to $A$, $B$ and $\varphi$.
Moreover, every element of finite order in $F$ is conjugated with some element in $G$ of finite order.
\end{fact}

Notice that the free product $G*\mathbb Z$ can be treated as an HNN-extension of $G$ by putting $A=B=\{e\}$ and $\varphi= \id$.

\begin{definition}
With the previous notation, we say that a word $w=g_0\ldots g_n$, which represents an element of $F$, is reduced if it contains no subwords:
  $tt^{-1}$; $t^{-1}t$; $tat^{-1}$ for an element $a\in A$; $t^{-1}bt$ for $b\in B$; $g h$ for $g,h \in G$; $e$, unless $w = e$.

A word is cyclically reduced if every cyclic permutation of it is reduced. An element of $F$ is cyclically reduced if it has a representation as a cyclically reduced word. 
\end{definition}

%%%Krzys: Dorzucilem warunki  na temat elementu neutralnego, ktore zazwyczaj pisze sie w tej definicji. W ostatniej linijce mozna by  sobie podarowac odrzucanie elementu neutralnego (bo wynika to z naszej definicji slowa zredukowanego), ale dodalem to dla podkreslenia.
\begin{definition}
With the previous notation, let $S_A \ni e$ and $S_B \ni e$ be transversals of right cosets of the groups $A$ and $B$ in $G$ .
We say that a reduced word $w=g_0\ldots g_n$ is in normal form if for every subword of $w$ of the form $t^{-1}g$ for $g \in G$ we have that $g\in S_B \setminus \{ e \}$,
and for every subword of the form $tg$ for $g\in G$ we have that $g\in S_A \setminus \{ e \}$.
\end{definition}
\begin{fact}
With the previous notation, every element of $F$ has a unique representation as a word in normal form.
\end{fact}

%%%Krzys1: Dodalem ponizszy komentarrz.
We will sometimes write reduced words as $a_0t^{\epsilon_1}a_1t^{\epsilon_2}a_2 \dots t^{\epsilon_m}a_m$, where $\epsilon _i \in \{-1,1\}$ and $a_i \in G$. Since our definition of reduced words do not allow neutral elements as letters, whenever some of the $a_i's$ in the above word are neutral, this word is treated as the word in which all these neutral $a_i$'s have been removed.

\begin{lemma}\label{aabb}
If $a, b \in G*\mathbb{Z}$, $e \neq ab=(ba)^{-1}$ and $G$ does not contain an involution, then $ab \in G^{G*\mathbb{Z}}$.
\end{lemma}
{\em Proof.} Wlog $b$ is cyclically reduced. The above equality is equivalent to $b^2 = (a^{-1})^2$.
 If $a$ is cyclically reduced, then $a\in G$ (otherwise there are no reductions after concatenation of the  words $a^{-1}$ and $a^{-1}$, and so, comparing letters, we get $b=a^{-1}$, a contradiction), 
therefore also $b\in G$, and we are done.
 If $a$ is not cyclically reduced, then, since $a^{-2}$ is cyclically reduced (as $b^2$ is such), $a^{-2}$ has to be equal to $e$. 
It is a contradiction, because, by Fact \ref{pf}, $G*\mathbb{Z}$ does not contain an involution. \hfill $\square$

\section{Some properties of regular fields}\label{fields}

 When we talk about a regular field $K$, without loss of generality we assume that $K=\C$ instead of $K$ is definable in $\C$. We start from the following basic observation.

\begin{proposition}\label{proposition:basic}
Let $K$ be a regular field. Then the function $x \mapsto x^n$ is onto $K$ for every natural number $n>0$, and also the function $x \mapsto x^q-x$ is onto $K$ if $q>0$ is the characteristic of $K$. In particular, $K$ is radically closed.
\end{proposition}

\noindent
{\em Proof.} Let $p$ be the global generic type of $K$ and $g \models p |_\emptyset$. Let $f\colon K \to K$ be one of the functions considered in the proposition. Then $g$ is algebraic over $f(g)$, so $p |_\emptyset$ cannot imply $p|_{f(g)}$ (as $p$ is non-algebraic). Hence, by the definition of regularity, $f(g) \models p|_\emptyset$, and so $f[K]$ is generic. Since this is a multiplicative or an additive subgroup, and regular groups are connected, we get that $f[K]=K$. \hfill $\square$\\

Let $K$ be a regular field, and let $p \in S_1(K)$ be its global generic type. If $K$ is generically stable, then $(K,cl_p)$ is an infinite dimensional pregeometry with the homogeneity property that for any finite $A \subset K$, $K \setminus cl_p(A)$ forms an orbit under $Aut(K/A)$. Therefore, \cite[Theorem 1.13]{HLS} or an obvious adaptation of the proof of \cite[Theorem 1]{KTW} yields

\begin{proven main theorem}\label{theorem:gen stable}
Each generically stable regular field is algebraically closed. In particular, each generically stable minimal or quasi-minimal field is algebriacally closed.
\end{proven main theorem}

Theorem \ref{theorem:gen stable} together with  Corollary \ref{corollary:SOP}  give us the following conclusion. 

\begin{corollary}\label{wd}
If the theory of a regular field $K$ is simple (or just has NSOP), then $K$  is generically stable, so algebraically closed.
\end{corollary}

%We will show now few applications of theorem \ref{alg}. 

The next corollary follows from Theorem \ref{theorem:gen stable} and \cite[Corollary 5.1]{PT} which says that the global generic type of a quasi-minimal structure of cardinality greater than $\aleph_1$ is generically stable. Working in the context of groups, we observed \cite[Corollary 5.1]{PT} independently in 2009, and we give here our short proof (not using the dichotomy theorem). Recall that a corresponding, well-known fact for minimal structures says that each uncountable minimal structure is strongly minimal.

\begin{corollary}\label{wp}
If $K$ is a quasi-minimal field and $|K| > \aleph_1$, then $K$ is generically stable, so algebraically closed.
\end{corollary}

\noindent
{\em Proof.} Let $p \in S_1(K)$ be the generic type of $K$. Since over any set of parameters $A$ of cardinality  $\le \aleph_1$ there are at most $\aleph_1$  definable, countable sets,
$cl_p(A)\ne K$. 
Hence, we can choose in $K$ a Morley sequence $(b_\alpha)_{\alpha<\omega_2}$ in $p$ over $\emptyset$. 
Let $\phi$ be such that $\models \phi (b_1,\ldots,b_n)$. 
In order to show that every permutation of elements $b_i$ satisfies $\phi$, it is enough to show this for a transposition of adjacent elements.
Since we can treat a prefix of the Morley sequence as parameters, we can focus on the transposition $(1,2)$. 
By the indiscernibility of the sequence $(b_\alpha)_{\alpha<\omega_2}$, we have that 
$ \phi(x,b_{\omega_1 2},\ldots,b_{\omega_1 n})$ is satisfied by all elements $b_{<\omega_1}$. Thus,  quasi-minimality implies that it is satisfied by co-countably many elements,
in particular by the element $b_\alpha$ for some $\omega_1 2<\alpha<\omega_1 3$. 
Since $\models \phi(b_\alpha,b_{\omega_1 2},\ldots,b_{\omega_1 n})$, 
the indiscernibility of the sequence $(b_\alpha)_{\alpha<\omega_2}$ gives us $\models \phi (b_2,b_1,b_3,\ldots,b_n)$. \hfill $\square$\\

%%%%%%%%%%%%%%%%%%%%%%%
\begin{comment}

\begin{corollary}\label{wd}
If the theory of field $K$ is simple (or just has $NSOP$) then $K$  is generically stable, so algebraically closed.
\end{corollary}
{\em Proof:} Theorems 4 and 5 from \cite{Pi} imply that in quasi-minimal field either generic type is generically stable 
or there exists a definable partial ordering such that every Morley sequence in p is strictly increasing. 
Second option is not possible by $NSOP$. ~~$\square$\\

\end{comment}
%%%%%%%%%%%%%%%%%%%%%

The above results lead to a question whether every quasi-minimal [regular] field is generically stable.
Such a question for groups was formulated by Pillay.
Tanovi\'c and independently the second author found a counter-example, which was later modified to the context of fields. 
An example of a quasi-minimal, non generically stable field can be found in \cite[Example 5.1]{PT}; it is an algebraically closed field with some ordering. A classification of minimal almost liner (so non generically stable) groups was found in \cite{KTW}. At the end of this paper, we give our original example of a non generically stable minimal group; an interesting property is that it has an elementary extension which is quasi-minimal.

These examples show that Theorem \ref{theorem:gen stable} is too weak to deduce Conjectures \ref{conjecture:quasi-minimal fields} or \ref{conjecture:regular fields} in their full generality (it was strong enough in the last two corollaries, because we used some extra assumptions). An interesting question is what happens if we additionally assume NIP. The example given at the end of the paper has NIP, so Theorem \ref{theorem:gen stable} seems to be too weak even in the NIP context (however, we have not checked whether the only known example of a non generically stable field \cite[Example 5.1]{PT} has NIP).
%Moreover this field has $NIP$, so we cant applay theorem \ref{alg} (in such a way as it was done in corollaries \ref{wp} i \ref{wd})
%to prove that quasi-minimal field with $NIP$ is algebraically closed.

Below we prove a theorem about regular fields which may have applications in proofs of algebraic closedness of fields 
that satisfy some additional assumptions (like NIP), which we discuss after the proof.

%\begin{theorem}\label{oo}
%Let $L$ be a finite, proper extension of field $K$
%Let  $\C \succ K$ be a monster model of $Th(K)$.
%By $\hat p$ we denote the global heir of $p$. 
%Then $\hat p^{(n)}$ has unbounded orbit under action of multiplicative group of field  $L(\C)$, where  $L(\C)$ is the interpretation of field $L$ in $\C$

\begin{theorem}\label{oo}
Let $K$ be a regular field with the global generic type $p$. Assume that $K$ has a proper, finite extension $L$ of degree $n$ (so $L$ can be identified with $K^n$ with the coordinate-wise addition and definable multiplication).
Then $p^{(n)}$ has unbounded orbit under the action of the multiplicative group of $L$.
\end{theorem}

\noindent
{\em Proof.} Suppose for a contradiction that $p^{(n)}$ has bounded orbit.
We will show that $K$ is generically stable, so, by Theorem \ref{theorem:gen stable}, it is algebraically closed, which yields
a contradiction with the assumption that $L$ was a proper, algebraic extension.

By Proposition \ref{proposition:basic}, $K$ is perfect. So
$L=K(a)$ for some $a$, and let $$f=X^n-b_{n-1}X^{n-1}-\ldots-b_0$$ be the minimal polynomial for $a$ over $K$. 
By the boundedness of the orbit, there exists an 
%%%Krzys: Zmienilem unbounded na uncountable.
uncountable set $X \subseteq K$
such that for every $\alpha, \beta \in X$ one has $$(\beta a+1)  p^{(n)}=(\alpha a+1)  p^{(n)}.$$ 
We take two elements $\alpha$ and $\beta$ from $X$ which are algebraically independent over the coefficients of $f$.
We have
$$
(\alpha a+1)^{-1}(\beta a +1) p^{(n)}=  p^{(n)}.
$$

For each element $a \in L$ the function $x \mapsto a \cdot x$ is a linear map from $L$ into $L$ treated as a vector space over $K$.
Choose $(1,a,a^2,\ldots,a^{n-1})$ as a basis of the vector space $L$ over $K$; then the elements considered above can be  viewed as the following matrices  (put $m=n-1$):

$$
\alpha a +1 =\left({ \begin{matrix} 
1 & 0 & 0 & 0 & \ldots & \alpha b_0 \\ 
\alpha & 1 & 0  & 0 & \ldots &\alpha b_1 \\
0 & \alpha   & 1 & 0  & \ldots & \alpha b_2 \\
\vdots & \vdots & \vdots &\vdots &\ddots & \vdots \\
0 & 0 & 0 & 0 & \ldots & 1+ \alpha b_m

 \end{matrix}} \right)~~~~
\beta a +1 =\left({ \begin{matrix} 
1 & 0 & 0 & 0 & \ldots & \beta b_0 \\ 
\beta & 1 & 0  & 0 & \ldots &\beta b_1 \\
0 & \beta   & 1 & 0  & \ldots & \beta b_2 \\
\vdots & \vdots & \vdots &\vdots &\ddots & \vdots \\
0 & 0 & 0 & 0 & \ldots & 1+ \beta b_m

 \end{matrix}} \right).
$$

Let $q_j= \sum_{i=0}^j (-\alpha)^{j-i}b_i~$ and $~h=\frac{1}{1+ \alpha q_m}$. Then $(\alpha a +1)^{-1}$ equals
%%%Krzys: Dopisalem minusy w ostatniej kolumnie.
$$
 \left({ \begin{array} {rrrrc}
1+(- \alpha)^{m+1} q_0 h &( -\alpha)^ {m} q_0 h &(- \alpha) ^{m-1} q_0 h & \ldots &  (-\alpha) q_0 h \\
 -\alpha+(- \alpha)^{m+1} q_1 h &1+( -\alpha)^ {m} q_1 h &(- \alpha) ^{m-1} q_1 h & \ldots  & (-\alpha) q_1 h\\
(-\alpha)^2 + (- \alpha)^{m+1} q_2 h &-\alpha+( -\alpha)^ {m} q_2 h &1+(- \alpha) ^{m-1} q_2 h & \ldots  & (-\alpha) q_2 h\\

\vdots & \vdots  &\vdots &\ddots & \vdots \\
(- \alpha)^m h & (- \alpha)^{m-1}h & (- \alpha)^{m-2}h & \ldots & h
 \end{array}} \right).
$$

 Let $M=(\alpha a+1)^{-1}(\beta a +1)$.
 The algebraic independence of $\alpha$ and $\beta$ implies that the entries of $M$ are non-zero, in particular
$$
M_{m-1,m}=\sum_{i=0}^{m-1}\beta b_i ((-\alpha)^{(m-1-i)} + (-\alpha)^{(m+1-i)}q_{m-1} h ) - (1+\beta b_m)(\alpha q_{m-1}h) \ne 0.
$$
We choose a Morley sequence $(c_0,\ldots,c_m)$  in $p$ over an arbitrary, but fixed, 
 small set of parameters $B$, containing the entries of $M$. 
Using the fact that $M \cdot  p ^{(n)}= p^{(n)}$ and Corollary \ref{corollary:dod}, we get
$$
\begin{array}{l}
\tp(c_{m-1},c_m/B) =\tp( (Mc)_{m-1}, (Mc)_m/B)= \tp(\sum_{i=0}^{m}M_{m-1,i} c_i, \sum_{i=0}^{m}M_{m,i} c_i/B)=\\
 \tp(\sum_{i=0}^{m}M_{m-1,i} c_i, \sum_{i=0}^{m-1}(M_{m,i}- \frac{M_{m-1,i}M_{m,m}}{M_{m-1,m}}) c_i/B).
\end{array}
$$
Let $d=\sum_{i=0}^{m}M_{m-1,i} c_i$ and  $e= \sum_{i=0}^{m-1}(M_{m,i}- \frac{M_{m-1,i}M_{m,m}}{M_{m-1,m}}) c_i$.
Since $c_m$ is generic over $c_{<m} \cup B$ and $M_{m-1,m}\neq 0$, 
by Corollary \ref{corollary:dod}, we have that $d$ is generic over $c_{<m} \cup B$, so also over $\{e\}\cup B$.  
We conclude that
$$
\tp(c_{m-1},c_m/B)= \tp(d,e/B)=\tp(c_m,c_{m-1}/B).
$$

We have proved that each 2-element Morley sequence in $p$ over any small set of parameters $B$ (containing the entries of $M$) is totally indiscernible over $B$ (i.e., its type over $B$  is invariant under the transposition). This implies that $p$ is generically stable. \hfill $\square$\\

%In other words, type of Morley sequence with two elements is invariant under changing order of variables  $(\spadesuit)$. 
%To show that Morley sequence over parameters is indiscernible set over $\emptyset$ it is enough to show that 
%every type $p^{(n)}$ is invariant under all permutations of variables.
%Moreover it enough to show that for transposition of two adjoining elements,
%so we need to show that for $c$ such that $\tp(c_0,c_1,\ldots ,c_m,c_{m+1})=p^{(n)}$ holds $\tp(c_0,c_1,\ldots ,c_m,c_{m+1})=\tp(c_0,c_1,\ldots ,c_{m+1},c_m)$.
%But this equality holds by $(\spadesuit)$ if we treat $b_{<m}$ as parameters.
 %$\square$\\

As it was mentioned, we hope that Theorem \ref{oo} may turn out to be useful in proving that some regular fields are algebraically closed.
The idea is of course to assume that our regular field $K$ has a proper, finite extension $L$ of degree $n$ and prove somehow 
that the orbit  $L \cdot p^{(n)}$ must be bounded (where $p$ is the global generic type of $K$).
%Then by theorem \ref{oo} we have that  $\hat p$ is generically stable, so field $K$ is algebraically closed by theorem \ref{alg}, 
%which is a contradiction with the assumption.

Let us look at this in the NIP context.
A (global) Keisler measure is a finitely additive, probabilistic measure $\mu$ on definable subsets of a given monster model $\C$.
A type $q$ is $\mu$-generic if for every $\varphi(x) \in q$,   $\mu(\varphi(\C))>0$.
From \cite{HPP}, it follows that in theories with NIP, there are boundedly many global $\mu$-generic types.
Assume $K$ is a regular field satisfying NIP and having a proper, finite extension $L$ of degree $n$.
If we could find a Keisler measure on $L$ invariant under the non-zero multiplicative translations (or satisfying a weaker condition that the action of the multiplicative group of $L$ moves $\mu$-positive sets 
onto $\mu$-positive sets) and such that $p^{(n)}$ is $\mu$-generic, then the type  $p^{(n)}$ would have bounded orbit, and we would be done. It is easy to show (see \cite[Section 4]{K}) that there exists a (definable) Keisler measure on $L$ invariant under the additive and non-zero multiplicative translations, but it is not clear how to construct such a measure $\mu$ so that $p^{(n)}$ is $\mu$-generic. 
%We would be able then apply theorem \ref{oo} and deduce that quasi-minimal fields with $NIP$ are algebraically closed.  

Let us finish with a comment on positive characteristic. A proof Podewski's conjecture in positive characteristic can be divided into two steps. First, using the Frobenius automorphism and the fact that minimal fields are radically closed, one gets that, working  in the pure field structure, $\acl(\emptyset)=F^{\alg}$ (i.e., the model-theoretic algebraic closure of $\emptyset$ coincides with the field-theoretic algebraic closure of the prime field), and this part goes through in the quasi-minimal context. Next, using Tarski-Vaught test, one checks that $\acl(\emptyset)$ is an elementary substructure of the field in question, which completes the proof. However, this part does not go through for quasi-minimal fields, as a co-countable set need not to intersect the countable filed $\acl(\emptyset)=F^{\alg}$.

\section{Examples around quasi-minimal groups}

We start from noticing that if there exists a counter-example to Conjecture \ref{conjecture:1} or Question \ref{question:2}, then there is a counter-example with only one non-trivial conjugacy class. 

Notice that if $G$ is a quasi-minimal group and $H$ is its proper, definable  subgroup, then $|H|\le \aleph_0$, 
because every coset of $H$ has the same cardinality as $H$.

\begin{lemma}\label{lemma:one class}
If there exists a non-abelian, quasi-minimal [regular] group $G$, then there exists such a group with exactly one non-trivial conjugacy class and so torsion-free. 
\end{lemma}

\noindent
{\em Proof.} Consider first the quasi-minimal case (in fact, it would be enough to give a proof in the regular case, but since the proof in the quasi-minimal case is very elementary and intuitive, we decided to include it). Since $Z(G)$ is a proper subgroup, it is countable.
It is easy to check that $G/Z(G)$ is quasi-minimal. Moreover, $G/Z(G)$ is centerless. Indeed, if  $xZ(G) \in Z(G/Z(G))$, then $x^G\subseteq xZ(G)$. 
Since the set of cosets of the centralizer of $x$ is in bijection with the conjugacy class of $x$, we see that $|C(x)|>\aleph_0$, and therefore $x\in Z(G)$. 
We conclude that replacing $G$ by $G/Z(G)$, we can assume that $G$ has trivial center. 

Let $x_1,x_2 \in G\backslash \{e\}$. 
They have countable centralizers, so uncountable (equivalently co-countable) conjugacy classes.
Hence, $ x_1^G \cap x_2^G \ne \emptyset$, so  $x_1^G=x_2^G$. Therefore, every two non-trivial elements are conjugated.

A standard argument shows that each group with only one non-trivial conjugacy class is torsion-free. Of course, all non-trivial elements have the same order which is infinite or equal to some prime number $p$. Assume that  the second possibility holds. Clearly $p \ne 2$, as otherwise $G$ is abelian, a contradiction. 
Let $x \in G\backslash \{e\}$ and $y$ be such that $yxy^{-1}=x^2$.
Then, by Fermat's little theorem, we have
\begin {equation}
x=y^pxy^{-p}= x^{2^p}=\left( {x^{2^{p-1}}}\right)^2= x^2.
\end{equation}
Hence, $x=e$, a contradiction.

Now, consider the regular case. Wlog $G=\C$ is the monster model in which we are working. By \cite[Proposition 3.7]{K1}, all non-central elements are conjugated. It remains to show that $G/Z(G)$ is regular. For this, we will check that if $(\tp(g/G),x=x)$ is strongly regular, then so is $(\tp(gZ(G^*)/G),``x \in G/Z(G)'')$, where $G^*\succ G$ is a bigger monster model containing $g$. Put $p=\tp(g/G)$ and $q=\tp(gZ(G^*)/G)$. First, notice that $q$ is non-algebraic. Otherwise a formula defining a union of finitely many cosets of $Z(G)$ belongs to the generic type $p$, so $Z(G)$ is generic, and hence $Z(G)=G$ by the connectedness of $G$, a contradiction. Since $p$ is invariant over $\emptyset$, so is $q$. Now, consider any small set of parameters $A$ and any element $aZ(G^*) \not\models q |A$. Then $a \not\models p |A$, so $p| A \vdash p|A,a$, which implies that $q |A \vdash q |A,aZ(G^*)$.\hfill $\square$\\

By virtue of Lemma \ref{lemma:one class}, a natural approach to find a counter-example to Conjecture \ref{conjecture:1} is to analyze methods of constructing groups with only one non-trivial conjugacy class. A fundamental tool in such constructions are HNN-extensions. A classical construction goes as follows:  

\begin{equation}
\begin{array}{l}
G_0 = F_{\aleph_1}, \\ 
G_{i+1}= \langle G_i \cup \{T_{st}: s,t \in G_i\backslash \{e\}\} | T_{st}sT_{st}^{-1}=t\rangle ,\\
G= \bigcup_{i\in \omega} G_i .
\end{array}
\end{equation}
\begin{remark}
The group $G$ constructed above is uncountable, it has one non-trivial conjugacy class, but it is not quasi-minimal as the centralizers of all non-trivial elements are uncountable. 
\end{remark}
{\em Proof:} By a basic property of HNN-extensions, we know that $G_0<G$, so $G$ is uncountable.
Every two elements of $G$ have to belong to $G_i$ for some $i$, so they are conjugated in $G_{i+1}$.

Let $r,s,t \in G_0\backslash\{e\}$. Then $T_{st}T_{rs}T_{tr} \in C_{G_1}(t)$. 
Since $G_0$ is uncountable, there exist uncountably many pairs $(r,s)\in G_0 \times G_0$. Moreover, for distinct pairs $(r_1,s_1) \ne (r_2,s_2)$ the elements $T_{s_1t}T_{r_1s_1}T_{tr_1}$ and $T_{s_2t}T_{r_2s_2}T_{tr_2}$ are also distinct.
Therefore, $C_{G}(t)$ is uncountable. 
Hence, $G$ is not quasi-minimal.\hfill $\square$\\

We do not know whether the above group is regular. A possible way to prove that it is not regular could be to show that there is no global generic type, or that there are more than one global generic types. Notice, however, that this group is not stable, as there are no stable groups with a unique non-trivial conjugacy class.

A reason why the above construction produces too big centralizers (and, in consequence, the group is not quasi-minimal) are ``redundant conjugations''. 
We modify this construction below, producing an uncountable group with a unique non-trivial conjugacy class and such that all non-trivial elements have countable centralizers, but we leave as an open question  whether this group is quasi-minimal or at least regular.\\

Let $G_0$ be a countable group with one non-trivial conjugacy class
(such a group can be constructed by starting in the above construction from a countable free group).
From now on we fix an element $x \in G_0 \backslash \{e\}$. For $H\geq G_0$,  $S_H$ denotes the conjugacy class $x^H$. 
By recursion, we construct a sequence of countable, torsion-free  groups $(G_\alpha)_{\alpha<\omega_1}$ which will satisfy:
\bi
\item[(i)] $G_\alpha \lneq G_{\alpha+1}$,
\item[(ii)] $\bigcup_{\alpha<\omega_1} G_\alpha = \bigcup_{\alpha<\omega_1} S_{ G_\alpha}$,
\item[(iii)] if $y \in S_{G_\alpha}$, then $C_{G_\alpha}(y)$ =$C_{G_{\alpha+1}}(y)$,
\item[(iv)]if $e \neq y \in G_\alpha\backslash S_{G_\alpha}$, then
 $\exists z_y \in G_\alpha \exists n \in \mathbb{N} (z_y^n=y \wedge \forall w \in G_{\alpha}\forall m\in \mathbb{N}\backslash \{0\} (wy^mw^{-1} \in \langle z_y\rangle \Rightarrow w \in \langle z_y \rangle) )$.
\ei

Let $G=\bigcup_{\alpha<\aleph_1} G_\alpha $, $Z_\alpha=\{z_y:y\in G_\alpha\backslash (S_{G_\alpha} \cup \{e\})\}$. 
Condition (i) implies that $G$ is uncountable, (ii) that there exists only one non-trivial conjugacy class, and (iii) that the centralizers of non-trivial elements are countable. 
The last condition will allow to maintain the others during the construction.

Suppose $(G_\alpha)_{\alpha<\beta}$ have been constructed so that Conditions (i) and (iii) hold for all $\alpha<\beta$ for which $\alpha+1<\beta$, Condition (iv) holds for all $\alpha<\beta$, and all $G_\alpha$'s are countable and torsion-free. We describe how to construct $G_\beta$.\\[2mm] 
{\bf Case 1 (limit step)} $\beta \in \LIM$.\\
We put $G_\beta:= \bigcup_{\alpha<\beta}G_\alpha$.
Conditions  (i) and (iii) are clearly satisfied, since they were satisfied for every $\alpha<\beta$; $G_\beta$ is clearly countable and torsion-free. 
Condition (iv) will be satisfied, since during the construction we will guarantee that for any $y$ considered in (iv) at some point of the construction, $z_y$ once chosen will be good at all further steps of the construction at which $y$ is considered in (iv). This will also justify the notation $z_y$ not depending on $\alpha$. \\[2mm] 
{\bf Case 2 (odd step)} $\beta=\alpha +1=\delta+ 2n+1$ for some $\delta \in \LIM\cup \{ 0 \}$ and $ n \in \mathbb{N}$.\\ 
We define $$G_{\alpha+1}:=G_{\alpha}*\mathbb{Z},$$ and we will check that our conditions hold. By Fact \ref{pf} and the assumption that $G_\alpha$ is countable and  torsion-free, we get that $G_{\alpha +1}$ is countable and torsion-free.
\bi
\item[(i)] is clearly satisfied.
\item[(iii)] Let $ y \in S_{G_\alpha}$. Suppose for a contradiction that there is some element  $k \in
C_{G_{\alpha+1}}(y) \backslash C_{G_\alpha}(y)$ and let $a_0t^{\epsilon_1}a_1\ldots t^{\epsilon_m}a_m$ be its normal form.
Then we have that
$a_0t^{\epsilon_1}a_1\ldots t^{\epsilon_m}a_m y a_m^{-1} t^{-\epsilon_m}\ldots  t^{-\epsilon_1}a_0^{-1}=y$. 
If  $a_m y a_m^{-1} \neq e$, then the term on the left defines an element from outside  $G_\alpha$ (because $m \geq 1$ since $k$ is outside $G_\alpha$), a contradiction.
The equality $a_m y a_m^{-1}= e$ implies that $y=e$, which is also a contradiction.
\item[(iv)] Let $y \in G_{\alpha}\backslash \{e\}$ and $y \notin  S_{G_{\alpha+1}}$. 
We will show that $z_y$ witnessing Condition (iv) for $G_\alpha$ will be good for $G_{\alpha+1}$. 
It is enough to show that for every $k \in G_{\alpha+1}\backslash G_\alpha$  and $m \in \mathbb N \setminus \{ 0 \}$ we have $ky^mk^{-1}\notin \langle z_y \rangle$. 
Arguing as in (iii), we get that if $ky^mk^{-1}\in \langle z_y \rangle \subset G_\alpha$, then $y^m=e$, and so $y=e$, which contradicts the choice of $y$.

It follows from the above paragraph that (iv) holds for any element $y \in G_\alpha^{G_{\alpha+1}} \setminus (S_{G_{\alpha+1}}\cup\{e\})$.

%%%Krzys: Zmienilem troche ponizszy dowod.
Let $y \in G_{\alpha+1}\backslash G_\alpha^{G_{\alpha+1}}$. 
%From the previous paragraph we know that $(iv)$ holds for elements in $G_\alpha^{G_{\alpha+1}}$, so wlog we can assume that $y\notin G_\alpha^{G_{\alpha+1}}$. 
Since conjugation is a group automorphism, we can assume that $y$ is represented as a cyclically reduced word $t^{\epsilon_1}a_1\ldots t^{\epsilon_s}a_s$.
Let $z_y$ be the shortest period of the word $t^{\epsilon_1}a_1\ldots t^{\epsilon_s}a_s$; in particular, we get $z_y^j=y$ for some $j \in \mathbb N\setminus \{ 0 \}$. 
Let $k$ be such that $ky^mk^{-1}=z_y^l$ for some $m \in \mathbb N \setminus \{ 0 \}$ and $l \in \mathbb{Z}$. 
We choose $r\in \mathbb N$ so that $|k|,|l| \ll r$. 
By comparing the numbers of letters $t^{\pm 1}$ in the equation $ky^{mr}k^{-1}=z_y^{lr}$ (notice that $z_y$ is cyclically reduced, so either $k$ fully cancels with a prefix of $y^{mr}$ and there are no cancellations between $y^{mr}$ and $k^{-1}$, or the other way around), we get that $l= \pm mj$. In  both cases below, we consider a long enough prefix or suffix of both sides of the equality $ky^{rm}k^{-1}=z_y^{lr}$, depending on where there are no cancellations.
If $l=-mj$, there would exist a division $z_y=ab$ such that $ba=z_y^{-1}$. 
By Lemma \ref{aabb}, we would get $z_y\in G_\alpha^{G_{\alpha+1}}$, so $y \in G_\alpha^{G_{\alpha+1}}$, a contradiction.
 So it must be the case that $l=mj$. 
%Since $z_y$ is aperiodic,  so by taking a large enough prefix of two side of equation $ky^{rm}k^{-1}=z_y^{lr}$, there exist $i,u$ such that $kz_y^i=z_y^u$.
%Hence $k=z_y^{u-i}$.
From the equation  $ky^{rm}k^{-1}=z_y^{lr}$, we conclude that either $k=z_y^i$ for some $i \in \mathbb Z$ and we are done, or there exists a proper division $z_y=ab$ such that $z_y=ba$, which implies that $z_y$ is periodic, a contradiction with the choice of $z_y$.
\ei
{\bf Case 3 (even step)} $\beta=\alpha +1=\delta+ 2n$ for some $\delta \in \LIM \cup \{ 0 \}$ and  $n \in \mathbb{N}\setminus \{ 0 \}$.\\
We define $$G_{\alpha+1}:=\langle G_\alpha, t| txt^{-1}=z\rangle$$ for some $z \in Z_\alpha$ (if $Z_\alpha=\emptyset$, i.e., $S_{G_\alpha}=G_\alpha \setminus \{ e\}$, then we could define $G_\alpha$ as in Case 2; but, in fact, $Z_\alpha$ is non-empty after every ``odd'' step).
At the end of the construction, we will describe how to choose the elements $z$ at these ``even'' steps of the construction in order to satisfy Condition (ii) after the whole construction. Now, we prove the other conditions.

As in Case 2,  by Fact \ref{pf} and the assumption that $G_\alpha$ is countable and torsion-free, we get that $G_{\alpha +1}$ is countable and torsion-free.
\bi
\item[(i)] obviously holds.

\item[(iii)] Since all elements from $S_{G_\alpha}$ are conjugated by elements from $G_\alpha$, it is enough to check (iii) only for $x$.
Suppose for a contradiction that there exists an element $k\in C_{G_{\alpha+1}}(x) \backslash C_{G_\alpha}(x)$, and write it in reduced form as
 $a_0t^{\epsilon_1}a_1\ldots t^{\epsilon_m}a_m$.
Then $a_0t^{\epsilon_1}\ldots t^{\epsilon_m}a_m x a_m^{-1} t^{-\epsilon_m}\ldots  t^{-\epsilon_1}a_0^{-1}=x$. 
If $\epsilon_1=1$, then $x$ would be conjugated  in $G_\alpha$ with some integer power $z^s$ for $s \ne 0$. 
Since $x$ is conjugated with $x^2$ in $G_\alpha$, $z^s$ would be conjugated with $z^{2s}$, which violates Condition (iv) for $G_\alpha$. 
Hence $\epsilon_1=-1$. 
Therefore, $a_1t^{\epsilon_2}\ldots t^{\epsilon_m}a_m x a_m^{-1} t^{-\epsilon_m}\ldots  t^{-\epsilon_2}a_1^{-1}=z^s$ for some non-zero $s \in \mathbb{Z}$. 
Analogously, we get that $\epsilon_2=1$, as otherwise $z^s$ would be conjugated in $G_\alpha$ with some power of $x$,
 and this would imply that $z^s$ and $z^{2s}$ are conjugated in $G_\alpha$. 
We conclude that $a_1$ conjugates two powers of $z$, so by (iv) for $G_\alpha$, $a_1$ must also be a power of $z$.
 It is a contradiction, since $t^{\epsilon_1}a_1t^{\epsilon_2}=t^{-1}a_1t$ can be reduced, which is impossible by the irreducibility of the word $a_0t^{\epsilon_1}a_1\ldots t^{\epsilon_m}a_m$.

\item[(iv)] Let $ y \in G_{\alpha}\backslash \{e\}$ and $ y \notin  S_{G_{\alpha+1}}$. 

We will show that $z_y$ from Condition (iv) for $G_\alpha$ will be good for $G_{\alpha+1}$. 
It is enough to show that for every $k \in G_{\alpha+1}\backslash G_\alpha$ and $m \in \mathbb N \setminus \{0\}$, if $ky^mk^{-1}\in \langle z_y \rangle$, then $y \in S_{G_{\alpha+1}}$. 
%Notice that the only way to reduce $t$ in word $ky^mk^{-1}$ is an application of relation $txt^{-1}=z$, and it implies that 
Since all letters $t^{\pm 1}$ in the word $ky^mk^{-1}$ must be reduced, 
$y^m$ is conjugated with a power of $x$ or of $z$ in $G_{\alpha}$.
It cannot be a power of $x$, because otherwise, since $x$ is conjugated with $x^2$ in $G_\alpha$, 
$y^m$  would be conjugated with $y^{2m}$, which would violate (iv) for $G_\alpha$. 
Hence, there exists $c\in G_\alpha$ such that $cy^mc^{-1}=z^l$ for some $l \in \mathbb Z \setminus \{ 0 \}$. 
Since $cyc^{-1}\in C(z^l)$ and it is easy to see that $z_z=z^{\pm 1}$ ($z_z$ chosen in (iv) for $G_\alpha$), by (iv) for $G_\alpha$, we get $cyc^{-1} \in \langle z\rangle$, in particular $y \in S_{G_{\alpha+1}}$, 
which contradicts our assumption.

It follows from the above paragraph that (iv) holds for any element $y \in G_\alpha^{G_{\alpha+1}} \setminus (S_{G_{\alpha+1}}\cup\{e\})$.
The bunch of lemmas below will finally allow us to deal with the remaining case when $y \in G_{\alpha+1}\backslash G_\alpha^{G_{\alpha+1}}$.

%%%Krzys; Ponizej dodatkowy lemat.
\begin{lemma}\label{dodatkowy}
If $a \in G_{\alpha+1} \setminus G_\alpha^{G_{\alpha+1}}$ and $n \in \mathbb N \setminus \{ 0 \}$, then $a^n \notin G_\alpha^{G_{\alpha+1}}$.
\end{lemma}
{\em Proof.} Suppose $a^n \in G_\alpha^{G_{\alpha+1}}$. Conjugating if necessary, we can assume that $a^n \in G_\alpha$. If $a^n \in S_{G_\alpha}$, then by (iii) and the fact that $a \in C_{G_{\alpha+1}}(a^n)$, we would get $a \in G_\alpha$, a contradiction. So $e \ne a^n \in G_\alpha \setminus S_{G_\alpha}$. Since $aa^na^{-1}=a^n \in G_\alpha$, the argument used above (in the first paragraph of the proof of (iv)) shows that $a^n \in S_{G_{\alpha +1}}$. Then, since $a \in C_{G_{\alpha+1}}(a^n)$, we see that $a$ is conjugated with an element centralizing $x$ which must lie in $G_\alpha$ by (iii), so $a \in G_\alpha^{G_{\alpha+1}}$, a contradiction. \hfill $\square$

\begin{lemma}\label{cent}
If $w,c\in G_{\alpha+1}$, $wc=cw$ and there exists $a\notin G_\alpha^{G_{\alpha+1}}$ commuting with $c$ and $w$, 
then there exists $r \in G_{\alpha+1}$ and $k,l \in \mathbb{Z}$ such that $w=r^k, c=r^l$.
\end{lemma}
{\em Proof.} By $|w|_t$ and $|c|_t$ we denote the total number of letters $t^{\pm 1}$ in reduced words representing $w$ and $c$, respectively.  
%The proof will be by induction on the total number of letters  $t^{\pm 1}$ in reduced words  representing  $w,c$.
The proof will be by induction on  $|w|_t +|c|_t$.

%%%Krzys: Odwolanie do dowodu powyzej.
{\em Base step:} Assume that $|w|_t=0$ or $|c|_t=0$, wlog $|c|_t=0$ . 
If $c\neq e$, then as in the proof of Lemma \ref{dodatkowy}, we get $c \in S_{G_{\alpha+1}}$ (otherwise $a$ would be in $G_\alpha$), so $a\in G_\alpha^{G_{\alpha+1}}$, a contradiction.
%Thus, we get that $a$ is conjugated with an element from the centralizer of $x$, and, by $(iii)$, this tells us that $a\in G_\alpha^{G_{\alpha+1}}$, a contradiction. 
We have proved that $c=e$, and we can put $r=w$.

{\em Induction step:} Wlog we can assume that $|w|_t \geq |c|_t$>0, $w$ is cyclically reduced, and $w,c$ are written in normal form. 
%%%Krzys: Zmienilem ponizsey uzasadnienie, a potem w ogole to usunolem, bo chyba z tego nie korzystamy.
%We can assume that the word $w$ ends with the letter $t$ (in other case we replace $w$ by $w^{-1}$).
%Conjugating if necessary, we can assume that the word $w$ ends with the letter $t$. 
%
%%%Krzys: Troche zmienilem reszte dowodu.
%%%%%%%%%%%%%%%%%%%%%%%%%%%%%%%%
\begin{comment}
Since $w$ is cyclically reduced,  after concatenation either of words $c,w$ or of $w, c^{-1}$ the only reduction which can be made 
is a multiplication of two elements of $G_\alpha$. 
If need we replace $c$ by $c^{-1}$, so we have such property for pairs $(w,c)$ and $(c,w)$.
 Wlog we can assume  that $|w|_t \geq |c|_t$. 
By uniqueness of normal form, equation $wc=cw$ provide us a division of word  $w= w_1w_2$  and $h\in G_\alpha$  such that  $w_2=hc$.
 Then $wc^{-1}$  commutes with $c$ and $a$, but $wc^{-1}=w_1 h^{-1}$, so it is shorter then $w$ in sense of letters $t^{\pm 1}$. From induction hypothesis we get  $wc^{-1}=r^k, c=r^l$, which implies $w=r^{l+l}$. $\square$
\end{comment}
%%%%%%%%%%%%%%%%%%%%%%%%%%%%%%%%
%
Since $w$ is cyclically reduced,  after concatenation either of words $c,w$ or of $w, c^{-1}$ the only possible reduction is a multiplication of two elements from $G_\alpha$. If there are no reductions of letters $t^{\pm 1}$ between $c$ and $w$, then, since $cw=wc$, the only possible reduction after concatenation of $w$ and $c$ can be a multiplication of two elements of $G_\alpha$. If there are no reductions of letters $t^{\pm 1}$ between $w$ and $c^{-1}$, then, since $c^{-1}w=wc^{-1}$, the only possible reduction after concatenation of $c^{-1}$ and $w$ can be a multiplication of two elements of $G_\alpha$. Replacing $c$ by $c^{-1}$ if necessary (and writing it in normal form), we have that $wc=cw$ and the only possible reductions after concatenation of $w, c$ and of $c,w$ can be multiplications of two elements of $G_\alpha$. By the uniqueness of normal form, the equation $wc=cw$ provides us a division $w_1w_2$ of the word $w$ and an element $h\in G_\alpha$  such that $w_2=hc$. Then $wc^{-1}$ commutes with $c$ and $a$, but $wc^{-1}=w_1 h$, so $|wc^{-1}|_t < |w|_t$. 
By induction hypothesis, we get  $wc^{-1}=r^k, c=r^l$, which implies $w=r^{k+l}$.\hfill  $\square$

\begin{lemma}\label{cykr}
Let  $\zeta \in G_{\alpha+1} \backslash G_\alpha^{G_{\alpha+1}}$, $a \in G_{\alpha+1}$, $n\in \mathbb{N}$ and $a=\zeta^n$. If $a$ is cyclically reduced, then so is $\zeta$. 
\end {lemma}
%%%Krzys: Poprawilem troche dowod tego lematu. Jest tu kwestia definicji postaci cyklicznie zredukowanej. Po pierwsze  element cyklicznie zredukowany to taki, ktory ma cyklicznie zredukowana reprezentacje; po drugie slowo, ktore zaczyna sie i konczy na cos z G_\alpha, nie jest cyklicznie zredukowane wedlug naszej definicji.
{\em Proof.} %Suppose for a contradiction that $\zeta$ is not cyclically reduced. Hence, 
We can write $\zeta$ as  $yby^{-1}$, where $b$ is cyclically reduced and such that we can choose a word representing $y$ and a cyclically reduced word $a_0t^{\epsilon_1}a_1\dots t^{\epsilon_m}$ representing $b$ so that the only possible reduction in the word $yby^{-1}$ can be a multiplication of two elements of $G_\alpha$ between $y, b$. Notice that $b \notin G_\alpha$, as otherwise $\zeta \in G_\alpha^{G_{\alpha +1}}$.
So, since $b$ is written as a cyclically reduced word,
the only possible reduction in $a= yb^ny^{-1}$ can be a multiplication of two elements from $G_\alpha$ between $y,b$.
Hence, if $|y|_t>0$, then $a$ is not cyclically reduced, a contradiction. So $y \in G_\alpha$. Consider the case $\epsilon _1=1$ and $\epsilon_m=1$ (the other cases are analogous). Since $b$ is cyclically reduced, the condition that $\zeta$ is not cyclically reduced as an element of the group (i.e., all its representations are non cyclically reduced words) is equivalent to the condition $ya_0 \notin \langle z \rangle$ and $y^{-1} \notin \langle x \rangle$. But this condition implies that $a=yb^{n}y^{-1}$ is not cyclically reduced (as an element of the group), a contradiction. Thus, $\zeta$ is cyclically reduced.\hfill $\square$

\begin{lemma}\label{ip}
If $a \in G_{\alpha+1} \backslash G_\alpha^{G_{\alpha+1}}$, then there exist $n\in \mathbb{N}$ and $w\in G_{\alpha+1}$ 
such that $w^n=a$ and $w$ does not have proper roots in $G_{\alpha+1}$. Such a $w$ we will call a minimal root.
\end {lemma}
{\em Proof.} In order to prove it, it is enough to find a bound on the degree of roots which we can take of $a$. 
Wlog $a$ is cyclically reduced, so, by Lemma \ref{cykr}, its every root $w$ is also cyclically reduced. 
Since $a \notin G_\alpha$, $|w|_t \geq1$.  
Hence, the degree of the root $w$ is bounded by $|a|_t$.\hfill $\square$

\begin{lemma}\label{nn}
If $a,b \in G_{\alpha+1} \backslash G_\alpha^{G_{\alpha+1}}$, $n\in \mathbb{N} \backslash \{0\}$ and $a^n=b^n$, then $a=b$.
\end{lemma}
{\em Proof.} Wlog we can assume that $a$ is cyclically reduced and ends with the letter $t$.
Then  $a^n$ is also cyclically reduced, so, by Lemma \ref{cykr}, $b$ is cyclically reduced.
By the uniqueness of normal form, the equation $a^n=b^n$ implies that there exists $p\in \mathbb{Z}$ such that $a=x^pb$.
Analogously, by the uniqueness of normal form for the term $a^{-n}=b^{-n}$, we have that there exists $m\in \mathbb{Z}$ such that $a=bx^m$.
We know that $x^p$ and $x^m$ are conjugated by some $c\in G_0$, i.e., $cx^mc^{-1}=x^p$. 
If $p\neq 0$, then, as $bc^{-1}\in C_{G_\alpha+1}(x^p)$, by (iii), we have $bc^{-1}\in G_0$.
Hence, $b \in G_0$, a contradiction. Thus, $p=0$, so $a=b$. \hfill $\square$

\begin{lemma}\label{jsc}
Assume that $a,b \in G_{\alpha+1} \backslash G_\alpha^{G_{\alpha+1}}$  and 
$n,m\in \mathbb{N} \backslash \{0\}$, and $a,b$ do not have proper roots in $G_{\alpha+1}$.
If $a^n=b^m$, then $a=b$ and $m=n$.
\end{lemma}
{\em Proof.} By Lemmas \ref{dodatkowy} and \ref{nn}, we can assume that $(n,m)=1$. 
Wlog $a$ is cyclically reduced and it ends with the letter $t$. Then, by Lemma \ref{cykr}, $b$ is also cyclically reduced.
The proof will be by induction on the number of letters $t^{\pm 1}$ in a shorter of the words $a,b$. Let $k=|a|_t$ and $l=|b|_t$.
 Since $a,b$ are cyclically reduced, $kn=lm$.

{\em Base step:} If $k=1$, then $n=lm$. Since $(n,m)=1$, we get $m=1$. As $b$ does not have a proper root, $n=1$, and we conclude that $a=b$.

{\em Induction step:} Wlog $k\leq l$. 
We can assume that $d:=NWD(k,l)< k$, because otherwise the equation $kn=ml$ implies $m|n$, so $m=n=1$ and $a=b$.
By Euclid's algorithm, there exist $p,q\in \mathbb{Z}$ such that $pk+ql=d$, $p>0$ and $q<0$. 
Since $a,b\in C_{G_{\alpha+1}}(a^n)$, also $a^pb^q\in  C_{G_{\alpha+1}}(a^n)$.  By Lemma \ref{dodatkowy}, $a^n \notin G_\alpha^{G_{\alpha+1}}$. Therefore,  by Lemma \ref{cent}, we get that $a^pb^q=r^{h_1}, a^n=r^{h_2}$ for some $r\in G_{\alpha+1}$ and $h_1,h_2 \in \mathbb Z$.  Moreover, we can assume that $r$ has no proper root by Lemma \ref{ip}.
 By Lemma \ref{cykr}, $r$ is cyclically reduced.
On the other hand, by the uniqueness of normal form and the equality $a^n=b^m$, we get that $b^{-q}$ is a suffix of $a^p$ modulo $x^i$ on the left side, and so $|a^pb^q|_t=pk+ql=d$.
%by $(iii)$ and $(iv)$ we can conclude that $a\notin G_\alpha^{G_{\alpha+1}}$, otherwise $a^n \in G_\alpha^{G_{\alpha+1}}$.
 By the above observations, $|r|_t\leq d <k$, so the equation  $a^n=r^{h_2}$  together with the induction hypothesis give us that $a=r^{\pm 1}$, and the equation $b^m=r^{h_2}$ that $b=r^{\pm 1}$.
 Since $r$ has infinite order, $a=b$ and $m=n$.\hfill $\square$

Now, we are ready to complete the proof of Condition (iv), i.e., we consider the case when $y\in G_{\alpha+1} \setminus G_\alpha^{G_{\alpha +1}}$.
%Let $y \in G_{\alpha+1}\backslash G_\alpha, y \notin  S_{G_\alpha+1}$. 
By Lemma \ref{ip}, there exists a minimal root $z_y$ of $y$, so $z_y^k=y$ for some $k\in \mathbb N$. 
Wlog we can assume that $y$ is cyclically reduced. By lemma \ref{cykr}, $z_y$ is also cyclically reduced. 
Let $c$ be such that $cy^mc^{-1}=z_y^n$ for some $m \in \mathbb N \setminus \{ 0 \}$ and $n \in \mathbb Z$. 
Raising this equation to a big power and counting occurrences of the letters $t^{\pm 1}$ on both sides, we easily get that $n=\pm mk$. 
Hence, for $d=c$ or $d=c^2$, we have $dy^m=y^md$. By Lemma \ref{dodatkowy} and  the assumption that $y \notin G_\alpha^{G_{\alpha+1}}$, we get $y^m \notin G_\alpha^{G_{\alpha+1}}$.
Thus,  Lemma \ref{cent} gives us that $d=r^{h_1}$ and $y^m=r^{h_2}$ for some $r \in G_{\alpha+1}$ and $h_1,h_2 \in \mathbb Z$. By Lemma \ref{ip}, we can assume that $r$ does not have proper roots.
 The equation $z_y^{mk}=r^{h_2}$ and Lemma \ref{jsc} imply that $r=z_y^{\pm 1}$. In the case when $d=c$, this shows that $c \in \langle z_y\rangle$. Consider the case $d=c^2$. Since $y^m=r^{h_2}$  and  $y^m \notin G_\alpha^{G_{\alpha+1}}$, we have $r\notin G_\alpha^{G_{\alpha+1}}$, hence, by Lemma \ref{dodatkowy}, $d \notin G_\alpha^{G_{\alpha+1}}$, and so $c\notin G_\alpha^{G_{\alpha+1}}$. Thus, $c$ has a minimal root $c_0$, say $c_0^l=c$. The equality $c_0^{2l}=r^{h_1}$ together with Lemma \ref{jsc} yield $c_0=r^{\pm 1}$, so $c \in \langle z_y \rangle$.
%If $d=c^2$, then analogically we get that $c=r^{\frac{h_1}{2}}$. Hence, $c\in \langle z_y \rangle$ 

\ei

Now, we will describe how to choose the elements $z$ at the ``even'' steps of the construction in order to satisfy Condition (ii). 
We always take this element $z\in Z_\alpha$ which was created as early as possible during the construction, but is still not conjugated with $x$.
Since all $G_\alpha$'s are countable, for each $\alpha<\omega_1$ after countably many steps all elements from $Z_\alpha$ will be conjugated with $x$. 
%It guarantee that all elements $z$ defined in condition $(iv)$ will be used. 
If $e \neq y \in G_\alpha \backslash x^{G_\alpha}$, then by (iv), there exists $z\in Z_\alpha$ and $ n\in \mathbb{N}$ such that $z^n=y$. 
But we know now that there exists $t\in G$ such that $txt^{-1}=z$. Hence $t x^nt^{-1}=y$. From the choice of $G_0$, we know that $x^n \in S_{G_0}=x^{G_0}$, so $y \in S_G$. 
Therefore, all non-trivial elements of the group $G$ are contained in the conjugacy class of $x$.\\

As it was mentioned in the previous section, we finish the paper with an example of a (non-pure) minimal group $G$ whose theory $T$ has quantifier elimination, there is a model $G_1\succ G$ which is quasi-minimal, and  the global generic type of $G$ (so also of $G_1$) is not generically stable; thus, a monster model of $\Th(G)$ is not generically stable. Moreover, $\Th(G)$ has NIP.  
%This gives us negative answers to \cite[Quastions 4.5, 4.7]{pi}. In fact, $G$ will be ordered of Type($\omega$) in the sense of Tanovi\'c, which shows that the existence of such an order in the context of groups does not lead to a contradiction with minimality.  

For the rest of the paper, $G=(F,0,+,<,P_n)_{n \in \omega}$, where:
\begin{itemize}
\item $(F,+,0)$ is the group of exponent 2 with neutral element $0$, spanned freely over $\mathbb{Z}_2$ by $(e_i: i \in \omega)$, i.e., $F=\bigoplus_{i \in \omega}\mathbb{Z}_2e_i$. In other words, $F$ is the subgroup of $\mathbb{Z}_2^\omega$ consisting of the elements with almost all coordinates equal to zero.
\item For $a \ne 0$, we define $$a<b \iff \max \{i \in \omega: \pi_i(a)=1\} < \max \{ i \in \omega : \pi_i(b)=1\},$$ where $\pi_i$ is the projection on the $i$-th coordinate. Moreover, $0<b \iff b \ne 0$.
\item For $n\ne 0$, we define $$P_n(a,b) \iff (a<b \wedge n=\max \{k \in \omega: \exists x_1,\dots,x_k(a<x_1<\dots<x_k<b)\}).$$ Moreover, $P_0(a,b)$ means that $a<b$ and there is no $c$ with $a<c<b$.
\end{itemize}
It is clear that all $P_n$'s are definable using $<$.

The proofs of the next lemmas are left as exercises.
\begin{lemma}\label{completeness}
The following list of axioms (which are all first order) axiomatizes $\Th(G)$.
\begin{enumerate}
\item $+$ is a group law with neutral element 0, giving a group of exponent 2.
\item $(\forall x \ne 0) (0<x)$.
\item The definitions of $P_n$'s.
\item The formula $\neg x<y \wedge \neg y<x$ defines an equivalence relation, which will be denoted by $\sim$.
\item $(\forall x,y,z) ((x\sim y \wedge x<z) \longrightarrow y<z)$ and $(\forall x,y,z) ((x \sim y \wedge z<x) \longrightarrow z<y)$. So, we can define $<^*$ on the quotient sort by $[x]_\sim <^*[y]_\sim \iff x<y$.
\item $<^*$ is a linear order with the smallest element $[0]_\sim$ and such that $[0]_\sim$ has an immediate successor, and any element different from $[0]_\sim$ has an immediate successor and an immediate predecessor.
\item $(\forall x,y) (x<y \longrightarrow x+y \sim y)$ and $(\forall x,y) (x \sim y \longrightarrow x+y<x)$. 
%Notice that this implies that for each $x\ne 0$, $Lin_{\mathbb{Z}_2}(y: y<x \lor y \sim x) = [x]_\sim \cup ([x]_\sim + [x]_\sim)$, and $dim_{\mathbb{Z}_2}([x]_\sim/Lin_{\mathbb{Z}_2}( y: y<x))=dim_{\mathbb{Z}_2}([x]_\sim/Lin_{\mathbb{Z}_2}(y: P_0(y,x)))=1$.
\end{enumerate}
\end{lemma}

\begin{lemma}
$\Th(G)$ has quantifier elimination.
\end{lemma}

The following corollary follows easily by quantifier elimination (to check NIP, one should use the fact that for theories with quantifier elimination it is enough to prove NIP for atomic formulas $\varphi(x,\overline{y})$, where $x$ is a single variable \cite[Proposition 12]{A}).

\begin{corollary}\label{G is minimal}
$G$ is minimal and has NIP.
\end{corollary}

The underlying order witnesses that a Morley sequence in the global generic type is not totally indiscernible, so the generic type is not generically stable.

Now, our goal is to construct a quasi-minimal model of $\Th(G)$. Notice that the order $(\omega_1,<)$ is not quasi-minimal, because the limit ordinals form a definable set which is neither countable nor co-countable. We will use a certain quasi-minimal modification of this order as a ``base'' of our structure. Namely, let $I$ be the order obtained from $(\omega_1,<)$ be replacing each infinite ordinal by a copy of $(\mathbb{Z},<)$. Let $G_1$ be the group of exponent 2 spanned freely over $\mathbb{Z}_2$ by $I$. For $a,b \in G_1 \setminus \{0\}$, we define $a<b \iff \max \{i \in I: \pi_i(a)=1\}<\max\{ i \in I: \pi_i(b)=1\}$, where $\pi_i$ is the projection on the $i$-th coordinate. Moreover, put $0<b \iff b \ne 0$. $P_n$'s are defined as in $G$. Then $G$ is a substructure of $G_1$.  It is clear that $G_1$ satisfies Axioms (1)-(7) from Lemma \ref{completeness}, so $G_1 \models \Th(G)$, and, by quantifier elimination, $G_1 \succ G$. By the choice of $I$,  a similar argument to the proof of Corollary \ref{G is minimal} yields the following

\begin{corollary}
$G_1$ is quasi-minimal.
\end{corollary}

\begin{question}
Does every minimal structure in a countable language posses a quasi-minimal elementary extension?
\end{question}

\vspace{10mm}
\noindent
{\bf Adresses:}\\
Tomasz Gogacz\\
Instytut Informatyki, Uniwersytet Wroc\l awski\\
ul. Joliot-Curie 15, 50-383 Wroc\l aw, Poland.\\
e-mail: gogo@cs.uni.wroc.pl\\[5mm]
Krzysztof Krupi\'nski\\
Instytut Matematyczny, Uniwersytet Wroc\l awski\\
pl. Grunwaldzki 2/4, 50-384 Wroc\l aw, Poland.\\
e-mail: kkrup@math.uni.wroc.pl

\end{document}